\documentclass[twoside,final,11pt]{amsart}

\usepackage{amsthm,latexsym,amssymb,amsfonts,amsmath,eucal,hyperref,enumerate}

\usepackage{color}           

\font\ninerm=cmr9
\font\nineit=cmti9
\font\ninebf=cmbx9

\newcommand{\textlineskip}{\baselineskip=13pt}
\newcommand{\smalllineskip}{\baselineskip=10pt}

\def\pmb#1{\setbox0=\hbox{#1}
    \kern-.025em\copy0\kern-\wd0
    \kern.05em\copy0\kern-\wd0
    \kern-.025em\raise.0433em\box0}

\def\ps@myheadings{\let\@mkboth\@gobbletwo
  \def\@oddhead{{\slshape\rightmark}\hfil{\footnotesize\thepage}}
  \def\@oddfoot{}
  \def\@evenhead{{\footnotesize\thepage}\hfil\slshape\leftmark}
  \def\@evenfoot{}
  \def\sectionmark##1{}\def\subsectionmark##1{}
}

\oddsidemargin=\evensidemargin
\def\runninghead#1#2{\pagestyle{myheadings}
\markboth{{\protect\footnotesize\it{\quad #1}}}
{{\protect\footnotesize\it{#2\quad}}}}
\newcommand{\fcaption}[1]{
        \refstepcounter{figure}
        \setbox\@tempboxa = \hbox{\footnotesize Fig.~\thefigure. #1}
        \ifdim \wd\@tempboxa > 5in
           {\begin{center}
        \parbox{5in}{\footnotesize\smalllineskip Fig.~\thefigure. #1}
            \end{center}}
        \else
             {\begin{center}
             {\footnotesize Fig.~\thefigure. #1}
              \end{center}}
        \fi}

\newcommand{\tcaption}[1]{
        \refstepcounter{table}
        \setbox\@tempboxa = \hbox{\footnotesize Table~\thetable. #1}
        \ifdim \wd\@tempboxa > 5in
           {\begin{center}
        \parbox{5in}{\footnotesize\smalllineskip Table~\thetable. #1}
            \end{center}}
        \else
             {\begin{center}
             {\footnotesize Table~\thetable. #1}
              \end{center}}
        \fi}

\newcommand{\bibit}{\nineit}
\newcommand{\bibbf}{\ninebf}
\renewenvironment{thebibliography}[1]
    {\frenchspacing
     \ninerm\baselineskip=11pt
     \begin{list}{\arabic{enumi}.}
        {\usecounter{enumi}\setlength{\parsep}{0pt}
    \setlength{\leftmargin 17pt}{\rightmargin 0pt}   
         \setlength{\itemsep}{0pt} \settowidth
    {\labelwidth}{#1.}\sloppy}}{\end{list}}

\numberwithin{equation}{section}
\newtheorem{theo}{Theorem}

\newtheorem{lemma}[theo]{Lemma}
\newtheorem{prop}[theo]{Proposition}
\newtheorem{cor}[theo]{Corollary}
\newtheorem{defi}[theo]{Definition}
\newtheorem{assum}[theo]{Assumption}

\theoremstyle{definition}

\newtheorem{rem}[theo]{Remark}

\numberwithin{theo}{section}

\numberwithin{equation}{section}

\newcommand\E{{\mathbb E}}

\renewcommand\Pr{\mathbb P}
\renewcommand\qed{\hbox{${\vcenter{\vbox{         
   \hrule height 0.4pt\hbox{\vrule width 0.4pt height 6pt
   \kern5pt\vrule width 0.4pt}\hrule height 0.4pt}}}$}}

\newcommand\R{\mathbb R}

\def\R{{\mathbb R}}
\def\E{{\mathbb E}}

\def\mv{{\rm v}}
\def\mv{{\rm e}}

\usepackage{eufrak}

\def\me{\mathsf{e}}
\def\mv{\mathsf{v}}

\begin{document}

\runninghead{S. Bonaccorsi, D. Mugnolo}{Long-time behavior of
  stochastically perturbed neuronal networks}

\normalsize\textlineskip
\thispagestyle{empty}
\setcounter{page}{1}

\centerline{\bf\Large Long-time behavior of stochastically perturbed}
\baselineskip=13pt
\par
\centerline{\bf\Large neuronal networks}
\baselineskip=13pt
\vspace*{0.37truein}

\centerline{\footnotesize Stefano BONACCORSI}
\baselineskip=12pt
\centerline{\footnotesize\it Dipartimento di Matematica, Universit\`a di Trento,}
\baselineskip=10pt
\centerline{\footnotesize\it Via Sommarive 14, 38050 Povo (Trento), Italia}
\baselineskip=10pt
\centerline{\footnotesize\it stefano.bonaccorsi@unitn.it}
\baselineskip=12pt
\centerline{\footnotesize Delio MUGNOLO}
\baselineskip=12pt
\centerline{\footnotesize\it Institut f\"ur Angewandte Analysis, Universit\"at Ulm, }
\baselineskip=10pt
\centerline{\footnotesize\it Helmholtzstra{\ss}e 18, D-89081 Ulm, Germany}
\baselineskip=10pt
\centerline{\footnotesize\it delio.mugnolo@uni-ulm.de}

\vspace*{0.225truein}

\vspace*{0.21truein}

\begin{abstracts} {Our investigation is specially motivated by the
    stochastic version of a common model of potential spread in a
    dendritic tree.  We do not assume the noise in the junction points
    to be Markovian. In fact, we allow for long-range dependence in
    time of the stochastic perturbation. This leads to an abstract
    formulation in terms of a stochastic diffusion with dynamic
    boundary conditions, featuring fractional Brownian motion. We
    prove results on existence, uniqueness and asymptotics of weak and
    strong solutions to such a stochastic differential equation. }
{diffusion on network, fractional Brownian motion, strong solutions of
  infinite dimensional stochastic differential equations, invariant
  measures.}
{34B45, 47D07, 60H15, 37L40}
\end{abstracts}

\vspace*{4pt}

\baselineskip=13pt

\normalsize

\section{Introduction}\label{sec:introduction}

This paper is concerned with a model of stochastically perturbed
parabolic network equations. We consider a network, i.e. a metric
graph, on whose nodes we impose continuity conditions, complemented by
dynamic or Kirchhoff-type laws that also incorporate stochastic
noise. We employ known results concerning well-posedness and
asymptotics of the deterministic model associated with the equations
in order to discuss solvability of the stochastic differential
equation and convergence properties of the stochastic convolution
process.

There is a well established theory concerning stochastic equations in
infinite dimensional spaces with additive Wiener noise, that can be
found for instance in \cite{dpz:Stochastic}. Our aim is to provide
similar results for perturbation driven by a fractional Brownian
motion. Since this noise is neither a semimartingale nor a Markov
process, a different approach to a stochastic calculus with respect to
it should be chosen.  Our first result, Theorem \ref{th:3.4}, provides
the existence of a weak solution for our problem. In this case we
extend the established results for Wiener noise. More interesting, and
less classical, is the existence of strong solutions provided in
Theorem \ref{th:3.6}. Even for the Wiener noise such a result, in the
infinite dimensional case, is known only in some special cases; as a
matter of fact, our assumption $H > 3/4$, together with the
ambientation $W(t) \in D((-A)^\alpha)$ for any $\alpha < 1/4$, which
leads to $H + \alpha > 1$, is the analog of the assumption $W(t) \in
D((-A)^\alpha)$ for some $\alpha > 1/2$ required in
\cite{barbu/daprato/roeckner/2007} (since the Wiener case corresponds
to $H = 1/2$). We shall also mention the paper \cite{KaLi08}, which
inspired our proof, but where the stronger requirement $W(t) \in D(A)$
is made.

\medskip

In Subsection~\ref{determ} we are going to introduce the deterministic
model of the time evolution of potential in a dendritical tree
according to Rall's linear cable theory. Successively, we are going to
perturb it by stochastic terms that are represented by a fractional
Brownian motion, as delineated in Subsection~\ref{stochas}. Section
\ref{sez3} is devoted to the stochastic analysis for a network which
contains only active nodes, while the general setting is sketched in
Section \ref{sez:4}.


\section{Biological motivations}
\label{sez2}

The motivation of this paper is provided by a stochastic model of
diffusion for neuronal activity in passive dendritic fibers. The
electrical behaviour of a neuron has been long studied, beginning with
the pioneering experiments of H. von Helmholtz, more than 150 years
ago. In the 1950s, the ground-breaking papers of, among other,
A.L. Hodgkin and A.F. Huxley~\cite{HH52} and W. Rall~\cite{Ra59} have
represented a major breakthrough in the modern neurophysiology. Ever
since, experimental results have dramatically increased. Successively,
more and more interest has been revived by stochastic models -- these
may help in gaining a better understanding of neuronal behaviour and
in predicting neuronal activity.

\medskip 

A typical neuron cell consists of four parts: the dendrites,
the synapses, the soma, and the axon. Dendrites are the input stage of
a neuron and receive input from other neurons through
\emph{axo-dendritic} synaptical junctions. Synapses also connect
dendrites to each other (leading to so-called \emph{dendro-dendritic}
interactions): in graph theoretical language, the dendritic network is
actually a tree converging into the soma, the cell's
body. In~\cite{Ra59, Ra60} Rall introduced the so-called ``lumped soma
model", representing (a part of) a dendritic tree as a lumped R-C
circuit, thus allowing for an internal electrical activity of the
soma. In Rall's model as well as in most of its generalizations
(and in particular in \emph{leaky integrate-and-fire neuronal models,
  see~\cite{Iz04})} it is the soma's duty to elaborate these inputs
and transmits them to other neurons through the axon.
In this way Rall was able, under strong (geometric) symmetry and
(electrotonic) isotropy assumptions, to describe the behaviour of a
whole dendritic tree by a cable equation on a single ``equivalent
cylinder" (i.e., an individual interval of finite length) equipped
with a dynamic boundary condition. Later on, multicylinder models have
been proposed already in the 1970s in~\cite{RR73}--\cite{RR74} in
order to weaken the theoretical assumptions on the topology of the
dendritic tree, thus describing more general classes of individual
neurons. Such models have been extensively discussed in the 1990s, see
e.g.~\cite{MEJ93} and subsequent papers. Aim of the present paper is
to further generalize a model of Rall-type to a non-deterministic one,
in order to allow for noise in the boundary, i.e., for
non-deterministic inputs in the somata.

\medskip

The starting point for the present analysis is represented by the
recent article~\cite{MR07}. There, a thorough investigation of the
deterministic version of such model has been performed: beside
well-posedness results, several qualitative properties of solutions
have been shown. Our mathematical approach also applies some of the
methods of that paper.

Rall's linear cable theory applies to the behaviour of the passive
dendritic tree: the propagation of action potentials along an axon is
modelled as a semilinear diffusion with a zero-order term that
accounts for describing how the cell expends energy in order to
propagate the signal, see below. On the other hand, it is known
that in axons action potentials are initiated and propagated with
fixed asymptotical speed and profile in the outward direction: an
effect of regenerative self-excitation is thus needed while modelling
axons.  Such a phenomenon is commonly described by a coupled system of
a semilinear diffusion equation and a nonlinear ordinary differential
equation. Such a system has been introduced by R. FitzHugh and J.
Nagumo as a simplification of the original Hodgkin--Huxley model: for
a survey of this theory we refer, e.g., to~\cite{Sc02}. Thus, the
behaviour of the neuron in the dendritic tree and the axon is
significantly different: as in~\cite{MEJ93} and~\cite{MR07} we will
not attempt to introduce axons into the model considered in this note.
A (naive) hybrid axo-dendritic model has been studied in a
deterministic setting in~\cite{CM07}.

Most of the aforementioned papers deal with a compartmental model of
the neuron, subsequently investigated by numerical methods. Aim of
this note is to provide a more theoretical approach based on an
interplay of operator theory and stochastic analysis. 

\medskip 

Stochastic fluctuations of synaptic activity and post-synaptic
elaboration of electronic potential is directly related to chaotic
interferences in the streaming of charge carrying molecules. This
chaotic streaming is classically modeled as a sequence of independent
random variables; averaging on time, in dependence of the time scale,
leads to choose a Brownian motion (or more generally a Levy process)
as a model for the stochastic input process. Possibly combined with
active properties of axons, this may lead to chaotic triggering of
action potentials. Several articles have actually dwelled on the
stochastic dynamics of axons within the framework of the
Hodgkin--Huxley or FitzHugh--Nagumo axonal models, see
e.g.~\cite{SD93} or~\cite{SFS98}, mainly concerned with the issue of
voltage fluctuations near threshold or the more theoretical approach
in \cite{BoMa-IDAQP}.

Less attention has seemingly been paid to the non-deterministic
aspects of sub-threshold stochastic behaviour, either in passive or
active fibers. However, there seem to be good reasons to perform such
an analysis: in particular, ``many computations putatively performed
in the dendritic tree (coincidence detection, multiplication, synaptic
integration and so on) occur in the sub-threshold regime" (quoted
from~\cite{MSK00}). A model for passive dendrites of infinite length
with distributed as well as synaptic current noise has been proposed
in \cite{kallianpur/wolpert, christensen/kallianpur}, whereas more
experimental studies have been presented in~\cite{DLK04}. A simplified
model of a network with active nodes is proposed in
\cite{BoMa-LDP}. Perhaps the clearest outcome of such investigations
is that \emph{voltage noise increases with depolarization},
cf.~\cite{JD+05}.


\medskip

Although our biological motivation leads us to consider a parabolic
network equation, 
our results (with possibly the exception of the existence of strong
solutions) are not confined to the 1-dimensional case. In fact, by
means of the same techniques we may consider with minor changes also a
stochastic initial-boundary value problem whose deterministic part is
a diffusion equation with dynamic (a.k.a. Wentzell-Robin) boundary
conditions on a bounded domain.  The same remark applies to the
stochastic analysis provided in Section \ref{sez3} which could have
been stated in abstract form and applied in different contests.


\subsection{Derivation of the equation}\label{determ}

The basic object of our investigations, i.e., the dendritic network is
modelled by a finite, connected graph $G$ with $m$ edges $\me_1,
\dots, \me_m$ and $n$ vertices $\mv_1,\dots,\mv_n$. 

Let $V$ denote the electric potential on the whole dendritic network,
including the junctions (synapses and further ramification points) and
the soma.  Then $V=V(t,x)$ is a function of the position $x$ along the
network and of time $t$ only. Of course, $V$ denotes the {\em
  deviation} from resting potential, which we rescale to 0 (in
concrete neurophysical measurements it is of approx.
$-70mV$). 
With a somewhat coarse but common approximation, we will assume the
cell's body to be isopotential.

By $V_j$ we denote the electric potential on the edge $\me_j$, i.e.,
$V_j(t,x) := V(t,\me_j(x))$.  Up to considering suitable rescaling
diffusion parameters in the equations, we may and do parametrize the
edges as intervals of unitary length. Actually, we avoid to introduce
more general network elliptic operators: in fact, under uniform
ellipticity assumptions this case is known to present no serious
mathematical challenges over the basic case of a plain network
Laplacian.

Following Rall's passive cable theory, the transmission of
post-synaptic potential in dendritic trees can be mathematically
described by a {linear cable equation}: on each edge $\me_j$ (i.e.,
each minimal set of dendritic elements that can be reduced to an
equivalent cylinder, cf.~\cite{MEJ93}), we consider the partial
differential equation
\begin{equation}\label{eq:bm1}
  \frac{\partial V_j}{\partial t}(t,x)
  = \frac{\partial^2}{\partial x^2}V_j (t,x) - p_j(x)V_j(t,x),\qquad t\geq 0,\; x\in (0,1),
\end{equation}

We borrow from graph theory some basic notions, in order to describe
the behaviour of the system at ramification points and to describe the
architecture of the dendritic network in a compact form. We consider a
graph ${\mathsf G}$ with vertex set ${\mathsf V}$ and edge set
${\mathsf E}$ throughout. We impose dynamical and time-independent
conditions in the nodes $\mv_1,\ldots,\mv_{n_0}$ and
$\mv_{n_0+1},\ldots, \mv_n$, respectively\footnote{In biological
  applications a dendritic network is in fact a tree, in the graph
  theoretical sense, and therefore only one soma needs to be
  considered, i.e., $n_0=1$.}. All nodes of the network are connected
by the above introduced edges $\me_j$.  We denote by ${\mathcal I} :=
{\mathcal I}^+ - {\mathcal I}^-$ the $n \times m$ \emph{incidence
  matrix} of the graph, where ${\mathcal I}^+:=(\iota^+_{ij})$ and
${\mathcal I}^-:=(\iota^-_{ij})$ are the \emph{incoming} and
\emph{outgoing} incidence matrices defined by
$$
\iota_{ij}^+ := \begin{cases}1& \text{if } \me_j(0) = \mv_i, \\ 0 &
  \text{otherwise,}\end{cases} \qquad \text{and} \qquad \iota_{ij}^-
:= \begin{cases}1& \text{if } \me_j(1) = \mv_i,\\ 0 &
  \text{otherwise,}\end{cases}
$$
respectively.  We also denote by $\Gamma(\mv_j)$ the set of all
indices of those edges having an endpoint in $\mv_j$, i.e.,
$$
\Gamma(\mv_i) = \{j \in \{1,\dots,m\}\ | \ \me_j(0)=\mv_i \text{ or
} \me_j(1) = \mv_i\}.
$$

Motivated by Rall's lumped soma model (see e.g.~\cite{Ra59}
and~\cite{MEJ93}), we specify the behaviour of the system at the
edges' boundaries by imposing several kinds of conditions in the
nodes.

\begin{itemize}
\item Both somata and synapses can be considered as isopotential. In
  other words, in all nodes $\mv_i$, $i=1,\ldots,n$, the potential $V$
  must satisfy the continuity assumption
  \begin{equation}\label{eq:bm3}
    V_j(t,\mv_i) = V_k(t,\mv_i)=:q^V_i \qquad \text{for all } j,k \in
    \Gamma(\mv_i),\qquad t\geq 0.
  \end{equation}

\item In the somata $\mv_1,\ldots,\mv_{n_0}$, the potential $q^V_i(t)
  = V(t,\mv_i)$ undergoes internal dynamics, subject to internal
  electrical activity and a stochastic feedback from the dendritic
  network. In the spirit of Rall's lumped soma model (see
  e.g.~\cite{MEJ93}) we impose general, possibly absorbing (and also
  possibly non-local) nodal conditions. They can be formulated as a
  Langevin-type equation
  \begin{equation}\label{eq:bm5}
    \frac{{\rm d}}{{\rm d}t} q^V_i(t) =  - \sum_{j=1}^m \iota_{ij} \frac{\partial
      V_j}{\partial x}(t,\mv_i) -\sum_{h=1}^{n}
    b_{ih}q^V_h(t)+\sum_{h=1}^n c_{ih} \dot Z_h(t), \qquad t\geq 0, 
  \end{equation}
  where the last addend on the right hand side is a sum, over the edges
  incident in $\mv_i$, of (formal) derivatives of the stochastic
  inputs acting there. We shall give a precise definition and complete
  assumptions on the noisy terms $Z^h$ in Section \ref{stochas}.

\item Motivated by the intrinsic randomness of external stimuli, we
  impose in the synapses and other ramification nodes $\mv_i$,
  $i=n_0+1,\ldots,n$, a version of a Kirchhoff's law which is also
  perturbed by a stochastic term $Z_i$. Such a condition can be
  equivalently formulated as
  \begin{equation}\label{eq:nonl}
    \sum_{j=1}^m \iota_{ij} \, \frac{\partial V_j}{\partial
      x}(t,\mv_i) + \sum_{h=1}^n b_{ih} q^V_h =\sum_{h=1}^n c_{ih}
    Z_h(t), \qquad t\geq 0.
  \end{equation}
\end{itemize}

Summing up, the biological model we want to discuss is governed by the
system of stochastic initial-boundary value problems
\begin{equation}
  \begin{cases}
    \dot{V}_j(t,x) = V_j''(t,x)-p_j(x) V_j(t,x), &t\geq 0,\;
    x\in(0,1), 
    \\
    &\phantom{t\geq 0,\;} j=1,\dots,m,
    \\
    V_j(t,\mv _i)= V_\ell (t,\mv _i)=: q^V_i, &t\geq 0,\; j,\ell\in
    \Gamma(\mv _i), 
    \\
    &\phantom{t\geq 0,\;} i=1,\ldots,n,
    \\
    \dot{q}_i^V(t) = -\sum\limits_{j=1}^m \iota_{ij} V'_j(t,\mv_i)
    -\sum\limits_{h=1}^{n} b_{ih} q^V_h + \sum\limits_{h=1}^n c_{ih}
    \dot{Z}_h(t), &t\geq 0,\; i=1,\ldots, n_0,
    \\
    \sum\limits_{h=1}^n b_{ih}q^V_h = -\sum\limits_{j=1}^m \iota_{ij}
    V'_j(t,\mv_i) + \sum\limits_{h=1}^n c_{ih} Z_h(t), &t\geq 0,\;
    i=n_0+1,\ldots,n,
    \\
    V_j(0,x) = V_{0j}(x), &x\in (0,1), \; j=1,\dots,m,
    \\
    V(0,\mv_i) = q_{0i}, &i=1,\ldots,n_0.
  \end{cases}
\end{equation}

For the sake of notational simplicity, let us introduce the
\emph{Kirchhoff operators} $K_a,K_p$ mapping $H^2(0,1;{\mathbb C}^m)$
into ${\mathbb C}^{n_0}$ and ${\mathbb C}^{n-n_0}$, respectively,
defined by
$$K_a V:=\begin{pmatrix}
-\sum\limits_{j=1}^m \iota_{1j} V'_j(t,\mv_1)\\
\vdots \\
-\sum\limits_{j=1}^m \iota_{n_0 j} V'_j(t,\mv_{n_0})
\end{pmatrix}\quad\hbox{and}\quad
K_p V:=\begin{pmatrix}
-\sum\limits_{j=1}^m \iota_{n_0+1\, j} V'_j(t,\mv_{n_0+1})\\
\vdots \\
-\sum\limits_{j=1}^m \iota_{n j} V'_j(t,\mv_{n})
\end{pmatrix}.
$$
Thus, the vector $K_a V$ (resp., $K_p V$) represents the differences
between incoming and outgoing flows in each of the active (resp.,
passive) nodes of the network.

We introduce the matrices
\begin{equation*}
  B_a = \Big( b_{ih} \Big)_{\substack{i=1,\dots,n_0 \\ h = 1,\dots,
      n}}, \quad B_p = \Big( b_{(n_0+i)h}
  \Big)_{\substack{i=1,\dots,n-n_0 \\ h = 1,\dots, n}} \quad
  \text{and} \quad B = \binom{B_a}{B_p} = \Big( b_{ih}
  \Big)_{\substack{i=1,\dots,n \\ h = 1,\dots, n}}
\end{equation*}
Roughly speaking, $B_a$ and $B_p$ encode the inhibitory and excitatory
properties of active and passive nodes under the influence of the
whole system. Also, in order to model the stochastic terms we
introduce the matrix
\begin{equation*}
  C_{aa} = \Big( c_{ih} \Big)_{\substack{i=1,\dots,n_0 \\ h = 1,\dots,
      n_0}}, 
\end{equation*}
as well as
\begin{equation*}
  C_a = \Big( c_{ih} \Big)_{\substack{i=1,\dots,n_0 \\ h =
      1,\dots, n}}, \quad C_p = \Big( c_{(n_0+i)h}
  \Big)_{\substack{i=1,\dots,n-n_0 \\ h = 1,\dots, n}} \quad \text{and}
  \quad C = \binom{C_a}{C_p} = \Big( c_{ih}
  \Big)_{\substack{i=1,\dots,n \\ h = 1,\dots, n}}
\end{equation*}
With no loss of generality, we may and do assume that the entries of
the matrix $C$ are non-negative real numbers and that there exist some
entries in $C_a$ which are strictly positive.  In particular, the
stochastic terms in \eqref{eq:bm5}--\eqref{eq:nonl} account for both
the intrinsic cellular noise and the external input signals. Observe
that although the possibility of non-local conditions does not seem to
be realistic in biological system, it may be interpreted as some form
of external boundary feedback control of the system -- e.g., with the
aim of stabilization.


\subsection{Fractional Brownian motion}\label{stochas}

As customary in network applications, we use stochastic processes in
order to model the input of our system. Altough the stochastic
processes used in this field were often assumed to be Markovian, other
considerations show that real inputs may exhibit long-range
dependence, i.e., the behaviour of the process at time $t$ does depend
on the whole history up to time $t$. Another property that stochastic
inputs show, also in telecommunication networks, is {\em
  self-similarity}: their behaviour is stochastically the same, up to
a space scaling, on changing the time scale. These features are common
in physiological systems and have been already applied to the
neuroscience, see e.g.~\cite{OTHB00},~\cite{CLCHH07}, and the
contributions in~\cite[Part II]{RanDin03}.  Among the processes which
exhibit these properties, we propose to model our stochastic input
with fractional Brownian motion (fBm for short). This class of
processes also verify {\em stationarity} of the increments and
continuity of trajectories.

A $1$-dimensional fractional Brownian motion is a centered Gaussian
process $\{B^H(t),\ t \ge 0\}$, such that
\begin{equation*}
  \E[B^H(t) B^H(s)] = \frac1{2H} (t^{2H} + s^{2H} - |t-s|^{2H}).
\end{equation*}
The constant $H$ is called {\em Hurst parameter} and  takes value in
$(0,1)$. The construction of fBm was proposed by Mandelbrot and van
Ness \cite{mandelbrot1968fbm} using the representation
\begin{equation*}
  B^H(t) = \int_0^t (t-s)^{H-1/2} \, {\rm d}B_s + \int_{-\infty}^0
  (t-s)^{H-1/2} - |s|^{H-1/2}) \, {\rm d}B_s
\end{equation*}
where $B = \{B_t,\ t \ge 0\}$ is a standard Brownian motion. We shall
be interested in a different characterization of fBm, based on the
observation that, given the fBm $B^H(t)$, there exists a unique
Brownian motion $W(t)$, adapted to the same filtration, and a kernel
$K_H(t,s)$ such that the identity 
\begin{equation*}
  B^H(t) = \int_0^t K_H(t,s) \, {\rm d}W(s).
\end{equation*}
holds. In case $H = \frac12$, $B^{1/2}(t) = W(t)$ is a standard
Brownian motion. In this paper we restrict our interest to a
particular class of fBm and impose the following.

\begin{assum}
  We assume that there exists a probability space $(\Omega,{\mathcal F},
  {\mathbb P})$ endowed with a filtration $\{{\mathcal F}_t,\ t \ge
  0\}$ that satisfies the standard assumptions, i.e., it contains all
  negligible subsets of $\Omega$ and it is right-continuous.

  On this space we are given a family of independent standard Brownian
  motions $\{W_h(t),\ t \ge 0\}_{h=1,\dots, n}$; then we consider, for
  each $h = 1,\dots, n$, the fBm
  \begin{equation*}
    Z_h(t) = \int_0^t K_H(t,s) \, {\rm d}W_h(s), \qquad t \ge 0.
  \end{equation*}
  We let $Z^a(t) = (Z_1(t), \dots, Z_{n_0}(t))$ and $Z^p(t) =
  (Z_{n_0+1}(t), \dots, Z_n(t))$ denote the stochastic inputs in
  active and passive nodes, respectively, and $Z(t) = (Z^a(t),
  Z^p(t))$ be the $n$-dimensional stochastic input process; we refer
  to $Z(t)$ as a $n$-dimensional fBm.

  The Hurst parameter $H$ belongs to $(\frac12,1)$.
\end{assum}


\subsection{Wiener integrals for fBm}

For any given Hilbert space $\Xi$, we are interested in the
introduction of Wiener-type integrals for $\Xi$-valued, square
integrable functions, defined with respect to a fractional Brownian
motion $Z(t)$. This is a known topic in literature and we follow the
stochastic calculus of variations approach, see the monograph
\cite{BiHuOkZh} and the references therein; we shall sketch below the
notation and the main results. 

In the Brownian case $H = \frac12$, the space of Gaussian random
variables defined as $\xi = \int_{\R} f(s) \, {\rm d}W_s$ is
isomorphic to the space $L^2(\R;\Xi)$. This isometry is
classically proved first for step functions and then is extended to
the whole space.
A similar construction was first investigated in \cite{PT:00} for the
fractional Brownian motion. Let ${\mathcal E}(\Xi)$ be the set of step
functions $\Phi: \R \to \Xi$
\begin{equation*}
  \Phi(t) = \sum_{i=0}^m \Phi_i \, {\bf 1}_{[t_i,t_{i+1})};
\end{equation*}
the stochastic integral $I(\Phi)$ is defined by setting
\begin{equation*}
  I(\Phi) = \int_{\R} \Phi(s) \, {\rm d}B^H(s) = \sum_{i=0}^m \Phi_i
  (B^H(t_{i+1}) - B^H(t_i))
\end{equation*}
and it is a Gaussian random variable with zero mean and covariance
that can be computed using the covariance matrix of $B^H$. The set
of Gaussian random variables defined by the elements of ${\mathcal
E}(\Xi)$ is a subset of the first Wiener chaos
\begin{equation*}
  {\mathcal H}_1 = \{X \in L^2(\Omega;\Xi) \,|\, \exists\, (f_n)_{n
    \ge 0}
  \subset {\mathcal E}(\Xi)\ :\ \lim_{n \to \infty} I(f_n) = X \quad \text{in
    $L^2(\Omega;\Xi)$} \}
\end{equation*}
generated by $B^H$.  The reproducing kernel Hilbert space (RKHS for
short) $\Lambda$ is defined as the closure of ${\mathcal E}(\Xi)$ with
respect to the scalar product
\begin{equation*}
  ({\bf 1}_{[0,t]}| {\bf 1}_{[0,s]})= R(t,s);
\end{equation*}
the mapping $I: {\bf 1}_{[0,t]} \to B^H(t)$ defines an isometry
between $\Lambda$ and ${\mathcal H}_1$.

As opposite to the Brownian motion case, for $H>\frac12$ the RKHS
$\Lambda$ for the fBm $B^H(t)$ cannot be identified with a space of
functions; however, $\Lambda$ contains (linear subspaces which are
isometric to) inner product spaces of functions.  Proceeding as in
\cite{PT:00} we can prove that this space can be identified with the
space $|{\mathcal H}|$ of functions in $L^1(\R;\Xi) \cap L^2(\R;\Xi)$
endowed with the norm
\begin{equation*}
  \| \phi \|^2_{|\mathcal H|} = \alpha_H  \int_{\R} \int_{\R}
  |t-r|^{2H-2} |\phi(r)|_{\Xi} |\phi(t)|_{\Xi} \, {\rm
    d}r \, {\rm d}t.
\end{equation*}
This space contains the class of step functions ${\mathcal E}(\Xi)$
and the isometry
\begin{equation*}
  \E \left| \int_{\R} \phi(t) \, {\rm d}B^H(t) \right|^2_{\Xi} = \| \phi
  \|^2_{|\mathcal H|}
\end{equation*}
holds.


\section{The abstract formulation}
\label{sez3}

Before discussing the complete stochastic differential equation in
Section \ref{sez:4}, we begin by considering the system with the
following simplified form of (\ref{eq:nonl}):
  \begin{equation}\label{eq:nonl-simplified}
    \sum_{j=1}^m \iota_{ij} \, \frac{\partial V_j}{\partial
      x}(t,\mv_i) + \sum_{h=1}^n b_{ih} q^V_h = 0, \qquad t\geq 0.
  \end{equation}
To begin with, we introduce the Hilbert spaces $X := (L^2(0,1))^m$
and $\partial X_a := {\mathbb C}^{n_0}$.  On the domain
\begin{equation*}
  D({\mathcal A}):=\left\{
    \begin{aligned}
      {\mathfrak v} := &\begin{pmatrix}V\\ q^V_a\end{pmatrix}\in
      (H^2(0,1))^m\times{\mathbb C}^{n_0} \quad \text{s.\ th.} \quad
      \exists\, q^V\in {\mathbb C}^n \text{ with } ({\mathcal
        I}^+)^\top q^V=V(0),
      \\
      &({\mathcal I}^-)^\top q^V=V(1),\
      (q^V_1,\ldots,q^V_{n_0})=(q^V_{a1},\ldots,q^V_{a{n_0}}), \text{
        and } K_p V=B_p q^V
    \end{aligned}
  \right\}
\end{equation*}
we define the operator $\mathcal A$ by
\begin{equation}\label{definA}
  {\mathcal A}{\mathfrak v} := 
  {\mathcal A}\begin{pmatrix}V\\ q^V_a\end{pmatrix} :=
  \begin{pmatrix}
    V''  - pV
    \\
    K_a V - B_a q^V 
  \end{pmatrix},
\end{equation}
where $V''$ denotes the vector $(V''_1,\ldots,V''_m)\in X$
consisting of (one-dimensional) second derivatives of the real-valued
entries of the vector-valued function $V=(V_1,\ldots,V_m)\in
(H^2(0,1))^m.$

One can directly check that the initial value problem associated with
\eqref{eq:bm1}--\eqref{eq:bm3}--\eqref{eq:bm5}--\eqref{eq:nonl-simplified}
can be equivalently formulated as an abstract stochastic Cauchy
problem
\begin{equation}
  \label{eq:0205-1}
  \left\{
    \begin{aligned}
      {\rm d}{\mathfrak v}(t) &= {\mathcal A}{\mathfrak v}(t) \, {\rm
        d}t + {\mathcal C}_a \, {\rm d}Z^a(t),\qquad t\geq 0,
      \\
      {\mathfrak v}(0) &= {\mathfrak v}_0,
    \end{aligned}
  \right.
\end{equation}
where the initial value is given by ${\mathfrak v}_0:=(V_0,q_0)^T\in
{\mathcal X} := X\times \partial X_a$ and ${\mathcal C}_a$ maps every
vector $q \in \partial X_a$ in $(0, C_{aa} q)^T \in {\mathcal X}$.

\vskip 1\baselineskip

Our first aim is to collect some results on the underlying
deterministic model, which will show well-posedness and further
qualitative properties of our system. Although formally new, the
results of this section can essentially be proved combining the ideas
of~\cite{Mu07} and~\cite{MR07}, where parabolic network equations with
passive only, non-local interactions and local, active node conditions
have been considered, respectively.

In order to prove the generation property of the operator ${\mathcal
  A}$, we introduce a Hilbert space
\begin{equation*}
  {\mathcal V} := \left\{
    \begin{aligned}
      {\mathfrak v} := &\begin{pmatrix}V\\ q^V_a\end{pmatrix}\in
      (H^1(0,1))^m\times{\mathbb C}^{n_0} \quad \text{s.\ th.} \quad
      \exists\, q^V\in {\mathbb C}^n \text{ with } ({\mathcal
        I}^+)^\top q^V=V(0),
      \\
      &({\mathcal I}^-)^\top q^V=V(1),\
      (q^V_1,\ldots,q^V_{n_0})=(q^V_{a1},\ldots,q^V_{a{n_0}})
    \end{aligned}
  \right\}
\end{equation*}
and the sesquilinear form ${\mathfrak a}: {\mathcal V} \times
{\mathcal V} \to {\mathbb C}$ defined by
\begin{equation*}
  {\mathfrak a}({\mathfrak u},{\mathfrak v}) := {\mathfrak
    a}_0({\mathfrak u},{\mathfrak v}) + {\mathfrak a}_1({\mathfrak
    u},{\mathfrak v}) := (U'|V')_{X} + \left((pU|V)_{X} + (B
    q^U|q^V)_{\partial X_a}\right).
\end{equation*}
We emphasize that this form is in general not symmetric, but it is
always densely defined. Mimicking the proofs of~\cite[Lemma~3.4]{Mu07}
and~\cite[Lemma~3.3]{MR07} one can prove that the operator associated
with $\mathfrak a$ is $\mathcal A$. Observe that $\mathfrak
a$ and hence $\mathcal A$ are self-adjoint if and only if the scalar
matrix $B$ is hermitian, i.e., if and only if mutual interaction of
the nodes are symmetric.

Combining the techniques presented in~\cite{MR07} and~\cite{Mu07} for
discussing dynamic and non-local boundary conditions, respectively, we
can prove the following.

\begin{prop}\label{poscontr}
  Let $B\in M_n({\mathbb C})$ and $p\in (L^1(0,1))^m$. Then the
  operator $\mathcal A$ generates a strongly continuous, analytic and
  compact semigroup $({\mathcal S}(t))_{t\geq 0}$ on the Hilbert space
  ${\mathcal X} := X \times \partial X_a$.
\end{prop}

\begin{proof}
  We first observe that the leading term in the form $\mathfrak a$,
  i.e., ${\mathfrak a}_0$, is clearly ${\mathcal X}$-elliptic and
  continuous. Furthermore, it follows from the Gagliardo--Nirenberg
  inequality that the space
  \begin{align*} 
    C({\mathsf G}) &:= \left\{
    \begin{aligned}
      {\mathfrak v} := &\begin{pmatrix}V\\ q^V_a\end{pmatrix}\in
      (C([0,1]))^m\times{\mathbb C}^{n_0} \quad \text{s.\ th.} \quad
      \exists\, q^V\in {\mathbb C}^n \text{ with } ({\mathcal
        I}^+)^\top q^V=V(0),
      \\
      &({\mathcal I}^-)^\top q^V=V(1),\
      (q^V_1,\ldots,q^V_{n_0})=(q^V_{a1},\ldots,q^V_{a{n_0}})
    \end{aligned}
  \right\}
  \\
  &\simeq \left\{ V\in (C([0,1]))^m \quad \text{s.\ th.} \quad
    \exists\, q^V\in {\mathbb C}^n \text{ with } ({\mathcal I}^+)^\top
    q^V=V(0),\ ({\mathcal I}^-)^\top q^V=V(1) \right\}
  \end{align*}
  of continuous functions over the network is embedded in an
  interpolation space of order $\frac{1}{2}$ between ${\mathcal V}$
  and ${\mathcal X}$, and obviously ${\mathfrak a}_1:C({\mathsf
    G})\times C({\mathsf G})\to {\mathbb C}$ is bounded. Thus, by a
  suitable perturbation argument, cf.~\cite[Lemma~2.1]{Mug08}, one
  concludes that also their sum $\mathfrak a$ is ${\mathcal
    X}$-elliptic and continuous. Thus, by~\cite[Prop.~1.51,
  Thm.~1.52]{Ou05} the associated operator $\mathcal A$ generates a
  strongly continuous, analytic semigroup. Moreover, by
  Rellich--Khondrakov's theorem the imbedding ${\mathcal
    V}\hookrightarrow {\mathcal X}$ is compact and accordingly the
  semigroup is compact.
\end{proof}

\medskip

\begin{rem}
  If coefficients $b_{ij}$ and $p_j(x)$ are real constant, resp.\ real
  valued functions, then ${\mathcal S}(t)$ maps real valued functions
  in real valued functions. thus, despite of the choice of complex
  valued function spaces, we may appeal to this remark in order to
  justify the application of our results to biological models.
\end{rem}

\medskip

All further properties of the semigroup essentially depend on the
matrix $B$: e.g., it is possible to characterize asymptotics of the
semigroup generated by $\mathcal A$ by means of positivity of the
$L^1$-coefficient $p$ and/or of the definiteness of the (in general
non-hermitian) matrix $B$, cf.~\cite[Thm.~2.3]{CMN08}.

\begin{prop}\label{stability}
  If $p\geq 0$ and $B$ is positive semidefinite, then the following
  assertions hold.
  \begin{enumerate}
  \item $({\mathcal S}(t))_{t\geq 0}$ is contractive, i.e.,
    $\|{\mathcal S}(t)\| \leq 1$ for all $t\geq 0$, and hence mean
    ergodic.
  \item If additionally $p=0$ and $B^* {\mathbf 1}\not=0$, or else if
    $B=0$ and $p\not\equiv 0$, then $({\mathcal S}(t))_{t\geq 0}$ is
    strongly stable, i.e., $\lim_{t\to\infty} {\mathcal
      S}(t){\mathfrak v}=0$ for all ${\mathfrak v}\in {\mathcal X}$.
  \item If additionally $p\geq p_0>0$ and/or $B$ is positive definite,
    then $(S(t))_{t\geq 0}$ is uniformly exponentially stable,
    i.e. $\|{\mathcal S}(t)\|\leq e^{-\omega t}$ for some $\omega >0$
    and all $t\geq 0$.
  \end{enumerate}
\end{prop}

By compactness of $({\mathcal S}(t))_{t\geq 0}$ and an abstract
criterion due to Ouhabaz, see~\cite[Thm. 2.6]{Ou05} we obtain the
following. We do not sketch here the proof and refer the reader
to~\cite{MR07} and~\cite{Mu07} for technical details.

\begin{cor}\label{co:3.4}
  If $p= 0$ and $B=0$, then $({\mathcal S}(t))_{t\geq 0}$ converges
  towards the (strictly positive) projection onto the $1$-dimensional
  eigenspace spanned by the constant vector ${\bf 1}\in
  X\times \partial X_a$.
\end{cor}

\medskip

Assume for a moment that ${\mathcal A}$ is negative definite. Then,
since ${\mathcal A}$ is a sectorial operator, we can introduce the
scale of interpolation spaces defined with respect to the domain of
the operator ${\mathcal A}$. The leading term of the form associated
with $\mathcal A$ is symmetric, i.e., $\mathcal A$ is self-adjoint up
to dropping the boundary conditions given by $B$, which are only
relevant for higher powers of $\mathcal A$.

More precisely, it follows by the spectral theorem that complex
interpolation spaces ${\mathcal X}_{\alpha} := [D(-{\mathcal
  A}),{\mathcal X}]_\alpha$ are equivalent to the domains of
fractionary powers $(-{\mathcal A})^{\alpha}$ of $-\mathcal A$.

Factoring $\mathcal A$ as in~\cite[Lemma~2.7]{Mug05} one can (up to
similarity) decouple its domain into the product $H_0(\mathsf
G)\times \partial X_a$. Here $H_0(\mathsf G)$ denotes the space of
$(H^2(0,1))^m$-functions with Dirichlet boundary conditions in the
active nodes $\mv_1,\ldots,\mv_{n_0}$.  For $\alpha<\frac14$, it
follows that $D((-{\mathcal A}^\alpha)$ is isomorphic to
$H^{2\alpha}(G) \times \partial X_a$.

In general, it is always possible to perform the same steps by
defining the complex interpolation spaces ${\mathcal X}_{\alpha} :=
[D(-\lambda I -{\mathcal A}),{\mathcal X}]_\alpha$ for a suitable
constant $\lambda$.


\subsection{Stochastic differential equations with additive fBm}
\label{sez:3.1}

In this section we solve the stochastic differential equation
\eqref{eq:0205-1}
\begin{equation}
  \left\{
    \begin{aligned}
      {\rm d}{\mathfrak v}(t) &= {\mathcal A}{\mathfrak v}(t) \, {\rm
        d}t + {\mathcal C}_a \, {\rm d}Z^a(t),\qquad t\geq 0,
      \\
      {\mathfrak v}(0) &= {\mathfrak v}_0,
    \end{aligned}
  \right.
\end{equation}
where  ${\mathfrak v}_0:=(V_0,q_0)^T\in
{\mathcal X} := X\times \partial X_a$.

\begin{defi}
Following the Wiener case treated for instance in [DPZ92] we call an
${\mathcal X}$-valued adapted process ${\mathfrak v} = \{{\mathfrak
  v}(t),\ t \in [0,T]\}$ a {\em strong solution} to (\ref{eq:0205-1})
if ${\mathfrak v}$ has a version such that
\begin{enumerate}
\item[(S1)] for almost all $t \in [0,T]$, $\Pr({\mathfrak v}(t) \in
  D({\mathcal A})) = 1$;
\item[(S2)] for any $t \in [0,T]$, $\displaystyle \Pr\left(\int_0^t
    |{\mathcal A}{\mathfrak v}(s)|^2_{{\mathcal X}} \, {\rm d}s < +\infty\right) =
  1$, and
\item[(S3)] for any $t \in [0,T]$ there holds
  \begin{equation}
    \label{eq:0205-7}
    {\mathfrak v}(t) = {\mathfrak v}_0 + \int_0^t {\mathcal
      A}{\mathfrak v}(s) \, {\rm d}s + {\mathcal C}_a \, Z^a(t).
  \end{equation}
\end{enumerate}
\end{defi}

\begin{defi}
Similarly, we define a {\em weak solution} to (\ref{eq:0205-1}) if
\begin{itemize}
\item[(W1)] for any $t \in [0,T]$, $\displaystyle \Pr\left(\int_0^t
    |{\mathfrak v}(s)|^2_{{\mathcal X}} \, {\rm d}s < +\infty\right) = 1$, and
\item[(W2)] for any $t \in [0,T]$ and any $y \in D({\mathcal A}^*)$ it
  holds
  \begin{equation}
    \label{eq:0205-2}
    \langle {\mathfrak v}(t),y \rangle = \langle {\mathfrak v}_0,y
    \rangle + \langle \int_0^t {\mathfrak v}(s) \, {\rm d}s, {\mathcal
      A}^*y \rangle + \langle {\mathcal C}_a \, Z^a(t), y \rangle.
  \end{equation}
\end{itemize}
\end{defi}


\subsection{Existence of weak solutions}

In order to provide the existence of a weak solution, our main
interest lies in the stochastic convolution process
\begin{equation*}
    W_{\mathcal A}(t) := \int_0^t {\mathcal
    S}(t-\sigma) {\mathcal C}_a \, {\rm d}Z^a(\sigma).
\end{equation*}
The key point is to find estimates for the variance, like
\begin{align*}
{\mathbb E}|W_{\mathcal A}(t)|^2_{\mathcal X^2}
&= \sum_{k=1}^{n_0} \left\| \left|{\mathcal S}(t-\sigma) {\mathcal
      C}_a e_k
\right|_{\mathcal X^2} \, {\pmb 1}_{(0,t)}(\sigma) \right\|^2_\Lambda
\end{align*}
where $\Lambda$ is the RKHS associated with the fBm $Z^a$. This
means that in case $H=1/2$ (the Wiener case) $\Lambda$ is the space
$L^2({\mathbb R})$ and in case $H > 1/2$ we can take instead of
$\Lambda$ the space $|{\mathcal H}|$.

Since ${\mathcal S}(t)$ is a semigroup,
there exist constants $M \ge 1$ and
$\omega \in \R$ such that $\|{\mathcal S}(t)\|_{L({\mathcal X})} \le
M e^{\omega t}$;
however, Proposition \ref{poscontr} implies that we can choose $M = 1$
and $\omega < 0$.
It follows that
\begin{align*}
  {\mathbb E}|W_{\mathcal A}(t)|^2_{\mathcal X^2} \le M \|{\mathcal
    C}_a\|_{HS} \left\| e^{\omega \sigma} \, {\pmb 1}_{(0,t)}(\sigma)
  \right\|^2_\Lambda;
\end{align*}
in the Wiener case, the norm is given by $\int_0^t e^{2 \omega s} \,
{\rm d}s = \frac{1}{2\omega}(e^{\omega t}-1) \le \frac{1}{2|\omega|}$;
on the other hand, in the fBm case, the estimate becomes
\begin{equation}
  \label{eq:stima-omega}
  \int_0^t \int_0^t e^{\omega s} e^{\omega r} |s-r|^{2H-2} \, {\rm d}s
  \, {\rm d}s \le \textit{Const.}\ \frac{1}{|\omega|^{2H}} |1 - e^{2\omega
    t}|
\end{equation}
and we get
\begin{equation}
  \label{eq:060508-1}
  \sup_{t \in [0,T]} {\mathbb E}|W_{\mathcal A}(t)|^2_{{\mathcal X}} \le
  \textit{Const.}\ \|{\mathcal C_a}\|^2_{HS}. 
\end{equation}

With a little effort, it is possible to be more precise on the
regularity (in space) of the stochastic convolution process. We remark
that the case $H=1/2$ is treated in \cite[Proposition
6.17]{dpz:Stochastic}.

\begin{lemma}
  Assume that the operator ${\mathcal A}$ is injective and
  dissipative.  Hence, we can consider the fractional powers of the
  operator ${\mathcal A}$ and the family of spaces ${\mathcal
    X}_\alpha$ with equivalent norm $\|v\|_{{\mathcal X}_\alpha} :=
  |(-{\mathcal A}^\alpha v|_{{\mathcal X}}$.

 Let $H > 1/2$ and assume that
  \begin{equation*}
    H - \alpha > 0.
  \end{equation*}
  Then the stochastic convolution process verifies
  \begin{equation}
    \label{eq:2307-1}
    \E \int_0^T \|W_{\mathcal A}(t)\|^2_{{\mathcal X}_\alpha} \,
    {\rm d}t \le \textit{Const.}
  \end{equation}
\end{lemma}

\begin{proof}
  We get
  \begin{align*}
    \E \int_0^T \|W_{\mathcal A}(t)|^2_{{\mathcal X}_\alpha} \, {\rm d}t \le
    \sum_{k=1}^{n_0} \int_0^T \left\| \left\|{\mathcal S}(t-\sigma)
        {\mathcal C}_a e_k \right\|_{{\mathcal X}_\alpha} \, {\pmb
        1}_{(0,t)}(\sigma) \right\|^2_\Lambda \, {\rm d}t
  \end{align*}
  By definition of the norm in ${\mathcal X}_\alpha$ we may compute
  \begin{align*}
    \left\| \left\|{\mathcal S}(t-\sigma) {\mathcal C}_a e_k
      \right\|_{{\mathcal X}_\alpha} \, {\pmb 1}_{(0,t)}(\sigma)
    \right\|^2_\Lambda &\le \int_0^t \int_0^t |(-{\mathcal A})^\alpha
    S(t-\sigma) {\mathcal C}_a e_k| \, |(-{\mathcal A})^\alpha
    S(t-\theta) {\mathcal C}_a e_k | \, |\sigma-\theta|^{2H-2} \, {\rm
      d}\sigma \, {\rm d}\theta
    \\
    &\le |{\mathcal C}_a e_k|^2 \int_0^t \int_0^t
    \frac{M_\alpha^2}{\sigma^\alpha \, \theta^\alpha}
    |\sigma-\theta|^{2H-2} \, {\rm d}\sigma \, {\rm d}\theta
    \\
    &\le M_\alpha^2 \, C_{\alpha,H}\ |{\mathcal C}_a e_k|^2
    t^{H-\alpha}
  \end{align*}
  and for any finite $T > 0$ the quantity in \eqref{eq:2307-1} is
  finite.
\end{proof}

\medskip

\begin{theo}\label{th:3.4}
  If ${\mathcal A}$ generates a $C_0$-semigroup $({\mathcal S}(t))_{t
    \ge 0}$, then there exists a unique weak solution to
  (\ref{eq:0205-1}) and it has the representation
  \begin{equation}
    \label{eq:0205-3}
    {\mathfrak v}(t) = {\mathcal S}(t){\mathfrak v}_0 + \int_0^t
    {\mathcal S}(t-s) {\mathcal C}_a \, {\rm d}Z^a(s), \qquad t \in [0,T].
  \end{equation}
\end{theo}

\begin{proof}
  Fix $t \in [0,T]$. Let us first notice the following identity which
  holds for every $\Phi \in C^1([0,t];D({\mathcal A}^*))$:
  \begin{multline}
    \label{eq:0205-4}
    \langle {\mathfrak v}(t), \Phi(t) \rangle - \langle {\mathfrak
      v}_0, \Phi(0) \rangle
    \\
    = \int_0^t \langle {\mathfrak v}(s), \dot \Phi(s) \rangle \, {\rm
      d}s + \int_0^t \langle {\mathfrak v}(s), {\mathcal A}^* \Phi(s)
    \rangle \, {\rm d}s + \int_0^t \langle \Phi(s), {\mathcal C}_a \,
    {\rm d}Z^a(s) \rangle.
  \end{multline}
  Assume now that $u$ is a weak solution to
  (\ref{eq:0205-1}). Choosing $\Phi(s) = {\mathcal S}^*(t-s) y {\pmb
    1}_{(0,t)}(s)$ in (\ref{eq:0205-4}) we get
  \begin{multline*}
    \langle {\mathfrak v}(t), y \rangle - \langle {\mathfrak v}_0,
    {\mathcal S}^*(t)y \rangle
    \\
    = \int_0^t \langle {\mathfrak v}(s), \tfrac{d}{d s} {\mathcal
      S}^*(t-s) y \rangle \, {\rm d}s + \int_0^t \langle {\mathfrak
      v}(s), {\mathcal A}^* {\mathcal S}^*(t-s) y \rangle \, {\rm d}s
    + \int_0^t \langle {\mathcal S}^*(t-s) y, {\mathcal C}_a\, {\rm
      d}Z^a(s) \rangle
  \end{multline*}
  that is
  \begin{equation*}
    \langle {\mathfrak v}(t), y \rangle - \langle {\mathcal
      S}(t){\mathfrak v}_0, y \rangle = \int_0^t \langle y, {\mathcal
      S}(t-s) {\mathcal C}_a\, {\rm d}Z^a(s) \rangle 
  \end{equation*}
  which implies uniqueness of the weak solution.

  We now proceed to show that the process $u$ defined in
  (\ref{eq:0205-3}) is a weak solution to (\ref{eq:0205-1}). Taking in
  mind formula (\ref{eq:0205-2}) we compute
  \begin{align*}
    \int_0^t \langle {\mathcal S}(\sigma){\mathfrak v}_0 & +
    \int_0^\sigma {\mathcal S}(\sigma - s) {\mathcal C} \, {\rm
      d}Z^a(s), {\mathcal A}^*y \rangle \, {\rm d}\sigma
    \\
    &= \int_0^t \langle {\mathfrak v}_0, {\mathcal S}(\sigma)
    {\mathcal A}^* y \rangle \, {\rm d}\sigma + \int_0^t \int_s^t
    \langle {\mathcal S}(\sigma-s) {\mathcal C}_a\, {\rm d}Z^a(s),
    {\mathcal A}^* y \rangle \, {\rm d}\sigma
    \\
    &= \int_0^t \tfrac{d}{d\sigma} \langle {\mathfrak v}_0, {\mathcal
      S}^*(\sigma) y \rangle \, {\rm d}\sigma + \int_0^t \langle
    \int_0^{t-s} \tfrac{d}{d\sigma} {\mathcal S}^*(\sigma) y \, {\rm
      d}\sigma \,,\, {\mathcal C}_a\, {\rm d}Z^a(s) \rangle
    \\
    &= \langle {\mathcal S}(t){\mathfrak v}_0, y \rangle - \langle
    {\mathfrak v}_0,y \rangle + \int_0^t \langle {\mathcal S}^*(t-s) y
    - y, {\mathcal C}_a\, {\rm d}Z^a(s) \rangle
    \\
    &= \langle {\mathcal S}(t){\mathfrak v}_0, y \rangle - \langle
    {\mathfrak v}_0,y \rangle + \int_0^t \langle y, {\mathcal S}(t-s)
    {\mathcal C}\, {\rm d}Z^a(s) \rangle - {\mathcal C}_a\, Z^a(t)
  \end{align*}
  and the representation (\ref{eq:0205-3}) yields to the result.
\end{proof}

\begin{cor}
  \label{co:0205-1}
  Assume that ${\mathcal A}$ is a {\em bounded} linear operator. Then the
  stochastic convolution process
  \begin{equation}
    \label{eq:0205-5}
    W_{\mathcal A}(t) := \int_0^t {\mathcal S}(t-s) {\mathcal C}_a\, {\rm d}Z^a(s)
  \end{equation}
  verifies
  \begin{equation}
    \label{eq:0205-6}
    W_{\mathcal A}(t) = \int_0^t {\mathcal A} W_{\mathcal A}(s) \,
    {\rm d}s + {\mathcal C}_a\, Z^a(t). 
  \end{equation}
\end{cor}


\subsection{Existence of a strong solution}


\begin{theo}
  \label{th:3.6}
  Let ${\mathcal A}$ be a sectorial operator and assume that
  ${\mathcal C}_a$ maps ${\mathcal X}$ into ${\mathcal X}_\alpha$ for
  any $\alpha < 1/4$.
  Assume that the Hurst parameter $H$ of the fBm $Z^a(t)$ verifies $H
  > 3/4$ and take 
  the initial condition ${\mathfrak v}_0 \in D({\mathcal A})$. Then
  there exists a unique strong solution to equation (\ref{eq:0205-1}).
\end{theo}

Up to rescaling, since we deal with a finite time interval $[0,T]$,
we can assume that  ${\mathcal A}$ is dissipative and injective.

\subsection*{Yosida approximations}

Before to proceed to prove the existence of a strong solution for
(\ref{eq:0205-1}), we need to introduce the family of operators
${\mathcal A}_n$, $n\in\mathbb N$, defined by
\begin{equation*} 
  {\mathcal A}_n = n {\mathcal A} R(n,{\mathcal A}) =
  n^2 R(n,{\mathcal A}) - n I.
\end{equation*}
This family has been introduced in Yosida's own proof of the
celebrated Hille--Yosida generation theorem. It enjoys remarkable
properties: the operators ${\mathcal A}_n$ are bounded; they commute
with one another as well as with $\mathcal A$ and with the resolvent
operators of $\mathcal A$; and finally
\begin{align*}
  \lim_{n \to \infty} {\mathcal S}_n(t)x = {\mathcal S}(t)x \quad
  \text{for all $x \in {\mathcal X}$ and} \quad \lim_{n \to \infty}
  {\mathcal A}_nx = {\mathcal A}x \quad \text{for all $x \in
    D({\mathcal A})$},
\end{align*}
by~\cite[Lemma~II.3.4]{EN00}.

For future references we need a lemma concerning the speed of
convergence of the semigroups generated from Yosida approximations to
${\mathcal S}(t)$.

\begin{lemma}
  \label{le:2307-1}
  Let ${\mathcal A}$ be a sectorial operator. Let also
  ${\mathcal C} \in L({\mathcal X}, {\mathcal X}_\alpha)$ for some
  $\alpha \in (0,1)$. Define, for any $n \in {\mathbb N}$, the
  sequence of functions $f_n(\sigma)$ setting
  \begin{equation*}
    f_n(\sigma) = \|{\mathcal A}_n ({\mathcal S}_n(\sigma) - {\mathcal
      S}(\sigma)) {\mathcal C}\|_{L({\mathcal X})}.
  \end{equation*}
  Then for given $T > 0$ it holds
  \begin{equation*}
    \lim_{n \to \infty} \sup_{\sigma \in (0,T)} |f_n(\sigma)| = 0.
  \end{equation*}
\end{lemma}

\begin{proof}
  Write
  \begin{align*}
    \|{\mathcal A}_n ({\mathcal S}_n(\sigma) - {\mathcal S}(\sigma))
    {\mathcal C}\|_{L({\mathcal X})} &= n \|{\mathcal A}^{1-\alpha}
    R(n,{\mathcal A}) ({\mathcal S}_n(\sigma) - {\mathcal S}(\sigma))
    {\mathcal A}^{\alpha} {\mathcal C}\|_{L({\mathcal X})}
    \\
    &\le n \|{\mathcal A}^{1-\alpha} R(n,{\mathcal A}) \|_{L({\mathcal
        X})} \|({\mathcal S}_n(\sigma) - {\mathcal
      S}(\sigma))\|_{L({\mathcal X})} \ \|{\mathcal A}^{\alpha}
    {\mathcal C}\|_{L({\mathcal X})}
  \end{align*}  
  We estimate, for $\omega = s({\mathcal A})$ the spectral bound of
  ${\mathcal A}$ and some $\epsilon > 0$, for any $n > \omega +
  \epsilon$,
  \begin{align*}
    \|{\mathcal A}^{1-\alpha} R(n,{\mathcal A}) \|_{L({\mathcal X})}
    &\le \int_0^\infty e^{-n t} \|{\mathcal S}(t) {\mathcal
      A}^{1-\alpha} \|_{L({\mathcal X})} \, {\rm d}t
    \\
    &\le \int_0^\infty e^{-n t} e^{(-\omega+\epsilon)t} t^{\alpha-1}
    \, {\rm d}t
    \\
    &\le \Gamma(\alpha) \frac{1}{(n - \omega + \epsilon)^{\alpha}}.
  \end{align*}
  We now recall that in \cite{Blunck} it is proved that
  $\|R(\lambda,A_n) - R(\lambda,A)\| \le \frac{(1+M_1)^2}{n}$. Define ${\mathcal
    S}(t)$ by the Dunford-Taylor integral
  \begin{equation*}
    {\mathcal S}(\sigma) = \frac{1}{2\pi \imath} \int_\Gamma e^{\lambda
      \sigma} R(\lambda,{\mathcal A}) \, {\rm d}\lambda.
  \end{equation*}
  Then we get
  \begin{align*}
    \|({\mathcal S}_n(\sigma) - {\mathcal S}(\sigma))\|_{L({\mathcal
        X})} &\le \left\| \frac{1}{2\pi \imath} \int_\Gamma
      e^{\lambda \sigma} [R(\lambda,{\mathcal A}) -
      R(\lambda,{\mathcal A}_n)] \, {\rm d}\lambda
    \right\|_{L({\mathcal X})}
  \end{align*}
  where $\Gamma$ is a piecewise smooth curve consisting of three
  pieces: a segment $\Gamma_1 = \{r e^{-\imath (\theta-\epsilon)}\,:\,
  1 \le r \le \infty\}$, a segment $\Gamma_3 = \{r e^{\imath
    (\theta-\epsilon)}\,:\, 1 \le r \le \infty\}$ and the arc
  $\Gamma_2 = \{e^{i \beta} \,:\, -(\theta-\epsilon) \le \beta \le
  (\theta-\epsilon)\}$. Take the norm inside the integral; it follows
    \begin{align*}
      \|({\mathcal S}_n(\sigma) - {\mathcal S}(\sigma))\|_{L({\mathcal
          X})} &\le \frac{1}{2\pi} \int_\Gamma |e^{\lambda \sigma}|
      \left\| [R(\lambda,{\mathcal A}) - R(\lambda,{\mathcal A}_n)]
      \right\|_{L({\mathcal X})} \, {\rm d}\lambda
      \\
      &\le \frac{(1+M_1)^2}{2 n \pi} \int_\Gamma |e^{\lambda \sigma}|
      \, {\rm d}\lambda = \frac{(1+M_1)^2}{2 n \pi} \int_{\Gamma_2}
      |e^{\lambda \sigma}| \, {\rm d}\lambda \le \frac{(1 +
        M_1)^2}{n}.
    \end{align*}
    
\end{proof}


\subsection*{Proof of Theorem \ref{th:3.6}}

It is possible to take the initial condition ${\mathfrak v}_0 = 0$,
possibly considering the process $\tilde {\mathfrak v}(t) = {\mathfrak
  v}(t) - {\mathcal S}(t) {\mathfrak v}_0$. Then we shall prove that
${\mathfrak v}(t) = W_{\mathcal A}(t)$ is a strong solution of
\eqref{eq:0205-1} with zero initial condition.

We notice that, by Corollary \ref{co:0205-1}, the identity
\begin{equation*}
  W_{\mathcal A_n}(t) = \int_0^t {\mathcal S}_n(t-s) {\mathcal C}_a\, {\rm
    d}Z^a(s) = \int_0^t {\mathcal A}_n W_{{\mathcal A}_n}(s) \, {\rm
    d}s + {\mathcal C}_a\, Z^a(t). 
\end{equation*}
holds. We shall now send $n$ to infinity in order to get (\ref{eq:0205-7}).

As a first step we prove that
\begin{equation}
  \label{eq:0205-8}
  \lim_{n\to \infty} \sup_{t \in [0,T]} \E |W_{{\mathcal A}_n}(t) -
  W_{\mathcal A}(t)|^2_{{\mathcal X}} = 0.
\end{equation}
Notice that
\begin{equation*}
  W_{{\mathcal A}_n}(t) - W_{\mathcal A}(t) = \int_0^t [{\mathcal
    S}_n(t-s) - {\mathcal S}(t-s)] {\mathcal C}_a \, {\rm d}Z^a(s)
\end{equation*}
and setting $q_k:={\mathcal C}e_k$ (where $\{e_k\}$ denotes the
canonical basis of $\partial X_a$) one concludes that
\begin{multline*}
  \E|W_{{\mathcal A}_n}(t) - W_{\mathcal A}(t)|^2_{\mathcal X} = \E
  \left| \sum_{k=1}^{n_0} \, \int_0^t [{\mathcal S}_n(t-s) - {\mathcal
      S}(t-s)] q_k \, {\rm d}Z_k(s) \right|^2_{\mathcal X}
  \\
  = \sum_{k=1}^{n_0} \| |[{\mathcal S}_n(\cdot) -
  {\mathcal S}(\cdot)] q_k|_{\mathcal X} {\pmb 1}_{(0,t)}(\cdot)
  \|^2_{|{\mathcal H}|}.
\end{multline*}
Hence, since it holds that $\psi(s) := [{\mathcal S}_n(s) - {\mathcal
  S}(s)] e_k|_{\mathcal X} $ converges to 0 as $n \to \infty$,
uniformly in $s\in[0,T]$ and
\begin{equation*}
  \|\psi(\cdot) {\pmb 1}_{(0,t)}(\cdot) \|^2_{|{\mathcal H}|} = \int_0^t \int_0^t
  |\psi(r)| |\psi(s)| |s-r|^{2H-2} \, {\rm d}s \, {\rm d}r \le
  \frac{T^{2H}}{H(2H-1)} \sup_{s \in [0,T]} |\psi(s)|^2,
\end{equation*}
the claim \eqref{eq:0205-8} holds.

\vskip 1\baselineskip

Next we show that $W_{\mathcal A}(t) \in D({\mathcal A})$,
$\Pr$-a.s., for all $t \in [0,T]$. In the case of a $Q$-Wiener
process $W(t)$, the thesis follows, for instance, from an
application of [DPZ92: Proposition 4.15]. In our case, the
assumptions can be modified: in fact, we require
\begin{equation}
  \label{eq:150508-2}
  \sum_{k=1}^{n_0}  \left\| |{\mathcal A}{\mathcal S}(\sigma) {\mathcal
      C}_a e_k|_{\mathcal X} \, {\pmb 1}_{(0,t)}(\sigma)
  \right\|^2_{|{\mathcal H}|} < +\infty.
\end{equation}
We can then mimic the proof of the aforementioned result by Da
Prato--Zabczyk and obtain the claimed result. Just
observe that
${\mathcal C}_a e_k$ belongs to
${\mathcal X}_{\alpha}$ for any $\alpha < \frac14$, hence
\begin{equation*}
  |{\mathcal A}{\mathcal S}(\sigma) {\mathcal C}_a e_k|_{\mathcal X} \le
  \textit{Const.}\ t^{\alpha-1}, \qquad t \in [0,T],\ \alpha <
  1/4.
\end{equation*}
A direct computation shows that
\begin{equation}
  \label{eq:060508-3}
  \int_0^T \int_0^T r^{\alpha-1} s^{\alpha-1} |r-s|^{2H-2} \, {\rm d}r \, {\rm d}s =
  \textit{Const.}\ T^{2\alpha + 2H - 2}
\end{equation}
provided $\alpha+H-1 > 0$; notice that the existence of an $\alpha$
which satisfies the bounds $1-H < \alpha < 1/4$
only holds if $H > 3/4$.

\vskip 1\baselineskip

Now, our last claim is to prove the following convergence
\begin{equation}
  \label{eq:0205-9}
  \lim_{n\to \infty} \sup_{t \in [0,T]} \E |{\mathcal A}_n
  W_{{\mathcal A}_n}(t) - {\mathcal A} W_{\mathcal A}(t)|^2_{\mathcal
    X} = 0.
\end{equation}
Together with \eqref{eq:0205-8}, this will imply that we can pass to
the limit in \eqref{eq:0205-6} and prove \eqref{eq:0205-7}, hence the
claim follows.

We split the estimate in two parts, since
\begin{equation}
  \label{eq:0305-1}
  \E |{\mathcal A}_n
  W_{{\mathcal A}_n}(t) - {\mathcal A} W_{\mathcal A}(t)|^2_{\mathcal
    X} \le 2 (\Phi_{n,1}(t) + \Phi_{n,2}(t))
\end{equation}
with
\begin{align*}
  \Phi_{n,1}(t) &= \E |{\mathcal A}_n W_{{\mathcal A}_n}(t) -
  {\mathcal A}_n W_{\mathcal A}(t)|^2_{\mathcal X}
  \\
  \Phi_{n,2}(t) &= \E |{\mathcal A}_n W_{{\mathcal A}}(t) - {\mathcal
    A} W_{\mathcal A}(t)|^2_{\mathcal X}.
\end{align*}
Notice that
\begin{align*}
  \E|{\mathcal A}_n W_{{\mathcal A}_n}(t) - {\mathcal A}_n W_{\mathcal
    A}(t)|^2_{\mathcal X} \le C , \left\| |{\mathcal A}_n
    \left({\mathcal S}_n(t-s) - {\mathcal S}(t-s) \right) \, {\mathcal
      C}_a |_{L({\mathcal X})} \right\|^2_{|\mathcal H|}
\end{align*}
and this term was treated in Lemma \ref{le:2307-1}; it holds that
$\sup_{s \in (0,t)} f_n(s) \to 0$ as $n \to \infty$ and  from Lebesgue
dominated convergence theorem it follows that $\Phi_{n,1}(t) \to 0$.

Next we proceed with the second term in (\ref{eq:0305-1}). We can
follow the same steps to get
\begin{align*}
  \Phi_{n,2}(t) = \E |{\mathcal A}_n W_{\mathcal A}(t) - {\mathcal A}
  W_{\mathcal A}(t)|^2_{\mathcal X} = \sum_{k=1}^{n_0} \left\| |[{\mathcal A} -
    {\mathcal A}_n] {\mathcal S}(t-s) \, {\mathcal C}_a e_k|_{\mathcal
      X} \right\|^2_{|\mathcal H|}
\end{align*}
Now since $[{\mathcal A} - {\mathcal A}_n]x \to 0$ as $n \to \infty$
for any $x \in D({\mathcal A})$, it is sufficient to show that each
term is uniformly bounded by an integrable function and apply
dominated convergence theorem in order to obtain the desired
convergence. Further, recall that ${\mathcal A}_n$ is a bounded
operator; hence, proceeding as above, we get the bound (compare with
\eqref{eq:060508-3})
\begin{equation*}
  \left\| |[{\mathcal A} - {\mathcal A}_n] {\mathcal S}(s) \,
    {\mathcal C}_a e_k|_{\mathcal X} {\pmb 1}_{(0,t)}(s)
    \right\|^2_{|\mathcal H|} < +\infty.\
\end{equation*}
In conclusion we obtain that both $\Phi_{n,1}(s)$ and $\Phi_{n,2}(s)$
converge to 0 as $n \to \infty$ uniformly in $s \in [0,T]$ and the
last claim is proved.


\subsection{Long time behaviour}
\label{sez4.1}

In this subsection we will be concerned with the asymptotic behavior of
$W_{\mathcal A}(t)$ for $t \to \infty$. Let us recall from
\cite[Theorem 11.7]{dpz:Stochastic} the equivalence between the
existence of an invariant measure and the boundedness of the
covariance operator $Q_t$ for $W_{\mathcal A}(t)$ at all times $t\geq
0$. In the truly dissipative case we obtain the following.

\begin{prop}\label{pr:dissip-active-nodes}
  Assume that $p \ge p_0 > 0$ and/or $B$ is positive definite. Then
  $W_{\mathcal A}(t)$ converges as $t \to \infty$ to a centered
  Gaussian random variable $\mu$. Furthermore, $\mu$ is an invariant
  measure for equation (\ref{eq:0205-1}).
\end{prop}

\begin{proof}
  The assertion follows if we prove that
  \begin{equation*}
    \sup_{t > 0} \mathop{Tr}Q_t
    < \infty,
  \end{equation*}
  where $Q_t$ is the covariance operator of $W_{\mathcal A}(t)$. We
  have discussed this quantity in previous section and we know that
  \begin{align*} 
    {\mathbb E}|W_{\mathcal A}(t)|^2_X &= \sum_{k=1}^{n_0} \left\| \left|
        {\mathcal S}(\sigma) {\mathcal C}_a e_k \right|_{\mathcal X}
      {\pmb 1}_{(0,t)}(\sigma) \right\|^2_\Lambda
  \end{align*}
  Now since in our assumptions, by Proposition~\ref{stability}.(3) the
  semigroup is uniformly exponentially stable, there exist $M \ge 1$
  and $\omega > 0$ such that
\begin{equation*}
  \mathop{Tr}Q_t \le M \, \|{\mathcal C}_a\|_{HS} \, \int_0^\infty
  e^{-2 \omega \sigma} \, {\rm d}\sigma \le \textit{Const.}
  \end{equation*}
  In case $Z_h$ is a fBm, we use in the space $\Lambda$ the norm
  $|{\mathcal H}|$ and we get (compare \eqref{eq:stima-omega}), for $t$ large,
  \begin{equation*}
    {\mathbb E}|W_{\mathcal A}(t)|^2_{\mathcal X} \le M
    \sum_{k=1}^{n_0} \left\| \left| {\mathcal S}(\sigma) {\mathcal
          C}_a e_k \right|_{\mathcal X} \right\|^2_{|{\mathcal H}|}
    \le \textit{Const.}\, \|C_a\|_{HS} \left[\frac{1}{\omega^{2H}} +
      \underbrace{\frac{t^{2H-2}}{\omega^2}}_{\to 0\ \text{as}\ t\to
        \infty}\right].
  \end{equation*}
This concludes the proof.
\end{proof}

We study now the uniqueness of the invariant measure and we make
precise the weak convergence of the associated transition
probabilities. Notice that we can appeal again to the theoretical
results in \cite[Ch.\ 11]{dpz:Stochastic}; hence, our results are
based on the properties of the semigroup ${\mathcal S}(t)$ obtained
before.

\begin{prop}\label{stoch-stability}
  Assume that $p\geq 0$ and $B$ is positive semidefinite; then the
  following assertions hold.
  \begin{enumerate}
  \item If additionally $p=0$ and $B^* {\mathbf 1}\not=0$, or else if
    $B=0$ and $p\not\equiv 0$, then there exists at most one invariant
    measure $\mu$ for (\ref{eq:0205-1}).
  \item If additionally $p\geq p_0>0$ and/or $B$ is positive definite,
    then there exists one invariant measure $\mu$ for
    (\ref{eq:0205-1}). Further, for all initial condition ${\mathfrak
      v}_0$ the solution ${\mathfrak v}(t)$ converges in law and in
    $L^2(\Omega;{\mathcal X}^2)$ to $\mu$.
  \end{enumerate}
\end{prop}

\vskip 1\baselineskip

In different terms, the situation that we obtain is that there exists
an equilibrium state (which can be thought of as resulting in the long
term behavior) for the neuronal network. It may be interesting, in
this connection, to study large deviations to estimate the probability
of onsets of chaotic impulses, compare \cite{ringach-malone:2007,
  BoMa-LDP}.

\vskip 1\baselineskip

The situation is different in case there is no dissipation acting in
the network. Then by Corollary \ref{co:3.4} the semigroup converges
toward a projection and for some $t_0$ there holds
\begin{equation*}
  \|{\mathcal S}(t) {\mathcal C}_a e_k\| \ge \frac12 |{\bf
    1}|_{\mathcal X}
  \sum_{k} | \langle {\bf 1}, {\mathcal C}_a e_k\rangle |  \qquad t > t_0
\end{equation*}
which is bounded below by the sum of the entries of $k$-th column of
the matrix $C_a$. Since at least one entry of $C_a$ is strictly
positive, as $t$ goes to infinity it follows
\begin{equation*}
  \sum_{k=1}^{n_0} \left\| \left| {\mathcal S}(\sigma) {\mathcal C}_a
      e_k\right| {\pmb 1}_{(0,t)}(\sigma) \right\|_{|{\mathcal H}|}
  \ge \textit{Const.}\ \left\| {\pmb 1}_{(t_0,t)}(\sigma)
  \right\|_{|{\mathcal H}|} \ge \textit{Const.}\ (t-t_0)^{2H}.
\end{equation*}
Therefore, the trace of the covariance operator $Q_t$ becomes
unbounded as $t \to \infty$ indipendently of the kind of noise that we
consider; this (negative) result is interesting enough to be recorded
in a proposition.

\begin{cor}\label{co:no-dissip-active-nodes}
  In case no dissipation occurs in the network ($p \equiv 0$ and $B
  \equiv 0$) then there does not exist an invariant probability
  measure for the system.
\end{cor}


\section{General stochastic perturbation}
\label{sez:4}

In this section we briefly discuss a possible approach to
inhomogeneous boundary value problem. This method has been proposed
in~\cite{KMN03} and has been extended to stochastic differential
equations in~\cite{BoZi}.

Let $\partial X^2_p := {\mathbb C}^{n-n_0}$ and define a
\emph{boundary operator} $R:D({\mathcal A})\to \partial X^2_p$ by
\begin{equation*}
  R{\mathfrak v}=R(V,q^V_a)^T:= K_p V-B_p q^V.
\end{equation*}
Consider the space
\begin{equation*}
  D({\mathbb A}):=\left\{
    \begin{aligned}
      v := &\begin{pmatrix}{\mathfrak v}\\ \phi^V \end{pmatrix}\in
      D({\mathcal A}) \times \partial X^2_p \quad \text{s.\ th.} \quad
      \exists\, q^V\in {\mathbb C}^n \text{ with } ({\mathcal
        I}^+)^\top q^V=V(0),
      \\
      &({\mathcal I}^-)^\top q^V=V(1),\
      (q^V_1,\ldots,q^V_{n_0})=(q^V_{a1},\ldots,q^V_{a{n_0}}), \text{
        and } R{\mathfrak v} = \phi^V
    \end{aligned}
  \right\}
\end{equation*}
and the operator on $X^2\times \partial X^2_a\times \partial X^2_p$
defined by
\begin{equation*}
  {\mathbb A}v :=
  \begin{pmatrix}{\mathcal A}{\mathfrak v} \\ 0\end{pmatrix}.
\end{equation*}
Then, it is known (cf.~\cite[\S~3]{KMN03}) that $\mathbb A$ generates
an analytic semigroup because $\mathcal A$ does so. Such a semigroup
is never stable, since it acts as the identity on the second component
of the vectors in $(X^2\times \partial X^2_a)\times \partial
X^2_p$. In fact, even if $({\mathcal S}(t))_{t\geq 0}$ is uniformly
exponentially stable, the semigroup generated by $\mathbb A$ is given
by
\begin{equation*}
  {\mathbb S}(t):=
  \begin{pmatrix}
    {\mathcal S}(t) & (I-{\mathcal S}(t)){\mathcal D}^{{\mathcal
        A},R}_0
    \\
    0 & I
  \end{pmatrix},
\end{equation*}
where ${\mathcal D}^{{\mathcal A},R}_0$ is a so-called abstract
Dirichlet operator, cf.~\cite[Lemma~3.3]{KMN03}.

We can now examine the case of a system featuring 
the presence of (stochastic) inputs in the passive nodes. The system
is described by the operator matrix ${\mathbb A}$ on the space
$(X^2\times \partial X^2_a)\times \partial X^2_p$. Let $Z(t)$ be the
$n$-dimensional stochastic process which models the input in the
(active and passive) nodes, and ${\mathbb C} = (0\ C_a\ C_p)^T$ be the
covariance operator of $Z(t)$, where $C_a$ and $C_p$ are the matrices
defined in Section 2.  Then the stochastic model can be written in the
form
\begin{equation}
  \label{eq:sist-completo}
  \begin{aligned}
    {\rm d}u(t) = {\mathbb A}u(t) \, {\rm d}t + {\mathbb C} \, {\rm
      d}Z(t)
    \\
    u(0) = u_0
  \end{aligned}
\end{equation}
which is solved in mild form by
\begin{equation}
  \label{eq:sol-sist-completo}
  u(t) = {\mathbb S}(t) u_0 + \int_0^t {\mathbb S}(t-s) {\mathbb C} \,
  {\rm d}Z(s).
\end{equation}
Proceeding as in Section \ref{sez:3.1} we obtain the following result.

\begin{theo}
  Assume that $u_0 \in (X^2\times \partial X^2_a)\times \partial
  X^2_p$. Then the process $u = \{u(t),\ t \in [0,T]\}$ is a weak
  solution of (\ref{eq:sist-completo}).

  Assume further that $u_0 \in D({\mathbb A})$ and $H > 3/4$. Then the
  weak solution is a strong solution of (\ref{eq:sist-completo}).
\end{theo}

We can be more precise in formula (\ref{eq:sol-sist-completo}) since
we know the representation of the semigroup ${\mathbb S}(t)$:
\begin{equation*}
  {\mathbb S}(t):=
  \begin{pmatrix}
    {\mathcal S}(t) & (I-{\mathcal S}(t)){\mathcal D}^{{\mathcal
        A},R}_0
    \\
    0 & I
  \end{pmatrix},
\end{equation*}
hence
\begin{align*}
  u(t) =
  \begin{pmatrix}
    u_a(t) \\ C_p Z(t)
  \end{pmatrix}
\end{align*}
where
\begin{align*}
  u_a(t) &= \int_0^t {\mathcal S}(t-s) \, {\mathcal C}_a \, {\rm
    d}Z(s) + \int_0^t (I-{\mathcal S}(t-s)){\mathcal D}^{{\mathcal
      A},R}_0 C_p \, {\rm d}Z(s).
\end{align*}
In this case, whenever the matrix $C_p$ is not identically zero, we
are again in presence of a gaussian process with finite trace class
covariance operator at any time but not bounded as $t \to \infty$ and,
similarly to Corollary \ref{co:no-dissip-active-nodes} we obtain the
following result.

\begin{cor}\label{co:passive-nodes}
  If the behaviour of passive nodes is affected by some stochastic
  inputs (i.e., the matrix $C_p$ is not identically zero), then there
  does not exist an invariant probability measure for the system.
\end{cor}

\end{document}